\documentclass[12pt]{article}
\usepackage{amsfonts}
\usepackage{latexsym}
\newtheorem{theorem}{Theorem}
\newtheorem{proposition}{Proposition}
\def\qed{\hfill$\Box$}
\begin{document}
\title{A Condition for Distinguishing Sceneries on Non-abelian Groups}
\author{Martin Hildebrand\footnote{Department of Mathematics and Statistics, University at Albany, State University of New York, Albany, NY 12222. E-mail: {\tt mhildebrand@albany.edu}}}
\maketitle
\begin{abstract}
A scenery $f$ on a finite group $G$ is a function from $G$ to $\{0,1\}$.
A random walk $v(t)$ on $G$ is said to be reconstructive if the
distributions of 2 sceneries evaluated on the random walk with uniform
initial distribution are identical only if one scenery is a shift of
the other scenery. Previous results gave a sufficient condition for 
reconstructivity on finite abelian groups. This paper gives a ready
generalization of this sufficient condition to one
for reconstructivity on
finite non-abelian groups but shows that no random walks on finite
non-abelian groups satisfy this sufficient condition.
\end{abstract}

\section{Introduction}

In \cite{fty}, Finucane, Tamuz, and Yaari considered the question of scenery
reconstruction on finite abelian groups and built upon results of
Matzinger and Lember~\cite{ml}. Finucane, Tamuz, and Yaari posed a number of
open questions. One question involves finding a sufficient condition for
reconstructivity for finite non-abelian groups similar to a condition proved
for finite abelian groups. In this paper, we shall develop such a condition,
but we shall also show that this condition is never satisfied if the group
is non-abelian. The techniques used involve Fourier transforms on the group
and the Plancherel formula.

Consider a graph with vertex set $V$ and edge set $E$. A function 
$f:V\rightarrow\{0,1\}$ is said to be a {\it scenery}. Let $v(t), t\in
{\mathbb N}$ be the position of a particle performaing a random walk on
this graph. If $f_1$ and $f_2$ are two sceneries, can an observation of
$\{f_1(v(t))\}$ and $\{f_2(v(t))\}$ distinguish which scenery was observed?
In this paper, we shall focus on the case where the vertices correspond to
elements of a finite group.

We shall define a random walk $v(t), t\in{\mathbb N}$ on a
finite group $G$ as follows. $v(1)$
has some distribution. Let $Z_1, Z_2, ...$ be i.i.d. random elements of $G$,
and let $v(t+1)=Z_tv(t)$ for $t\in{\mathbb N}$. The step distribution
$\gamma$ of this random walk is given by $\gamma(s)={\mathbb P}(Z_n=s)$.

{\it Definition:} Let $\gamma:G\rightarrow{\mathbb R}$ be the step 
distribution of a random walk $v(t)$ on a finite group $G$ so that $v(1)$
is picked uniformly from the elements of $G$. Note that
 $\gamma(k)={\mathbb P}(v(t+1)v(t)^{-1}=k)$.
Then $v(t)$ is said to be {\it reconstructive} if the distributions
of $\{f_1(v(t))\}_{t=1}^{\infty}$ and $\{f_2(v(t))\}_{t=1}^{\infty}$ where
$f_1$ and $f_2$ are sceneries on $G$ are
identical only if $f_1$ is a shift of $f_2$, i.e. there exists a
$g\in G$ such that $f_1(k)=f_2(kg)$ for all $k\in G$.

To generalize the sufficient condition in \cite{fty}, we use
representation theory of finite groups and Fourier analysis.
Let $\rho$ be a representation of $G$; in other words, $\rho$ is a 
function from $G$ to $GL_n({\mathbb C})$ for some positive
integer $n$ such that $\rho(st)=\rho(s)\rho(t)$ for all $s,t\in G$.
The value $n$ is called the degree of $\rho$ and is denoted $d_{\rho}$.
Define the Fourier transform of $f$ on $\rho$ by
$\hat f(\rho)=\sum_{s\in G}f(s)\rho(s)$.
A representation $\rho$ is said to be irreducible if the only subspaces
$H$ of ${\mathbb C}^n$ such that $\rho(s)H=H$ for all $s\in G$ are
the zero subspace and ${\mathbb C}^n$ where $n=d_{\rho}$. Representations
$\rho_1$ and $\rho_2$ are said to be equivalent if for some invertible
matrix $A$, $\rho_1(s)=A\rho_2(s)A^{-1}$ for all $s\in G$. For more
details on representation theory, see chapter 2 of Diaconis~\cite{di}
or Serre~\cite{se}.

The generalization of the sufficient condition in \cite{fty} (Theorem 1.2)
is given by the following.

\begin{theorem}
\label{condition}
Let $\gamma$ be the step distribution of a random walk $v(t)$ on a finite
group $G$. Then $v(t)$ is reconstructive if the following condition holds:

If $\sum d_{\rho_1}...d_{\rho_n} Tr((\hat\gamma(\rho_1)^{\ell_1}
\otimes ... \otimes\hat\gamma(\rho_n)^{\ell_n})\hat J_n(\rho_1\otimes
...\otimes\rho_n))=0$ for all positive integers $\ell_1,...,\ell_n$ 
implies that $J_n(x_1,...,x_n)=0$ for all $x_1,...,x_n\in G$ where the
sum is such that $\rho_1$, ..., $\rho_n$ each range over all irreducible
representations of $G$ up to equivalence, then $f$ can be reconstructed up
to shifts.
\end{theorem}

We also show
\begin{theorem}
\label{muchado}
If $G$ is a non-abelian group, then no random walk $v(t)$ satisfies the
condition of Theorem~\ref{condition}.
\end{theorem}

\section{Proof of Theorem~\ref{condition}}

The proof of this theorem is an adaptation of the proof of Theorem A.1 in
\cite{fty}. The proof in \cite{fty} uses Fourier analysis on finite abelian
groups; the proof here uses Fourier analysis on groups which may be
non-abelian.

Define the {\it spatial autocorrelation} $a_{\rho}(\ell)$ for $\ell\in G$ by
\[
a_f(\ell)=\sum_{k\in G}f(k)f(\ell k).
\]
Its Fourier transform is given by the following proposition.
\begin{proposition}
$\hat a_f(\rho)=\hat f(\rho)\hat g(\rho)$
where $g(s)=f(s^{-1})$.
\end{proposition}

{\it Proof:}
\begin{eqnarray*}
\hat a_f(\rho)&=&\sum_{\ell\in G}a_f(\ell)\rho(\ell)\\
&=&\sum_{\ell\in G}\sum_{k\in G}f(k)f(\ell k)\rho(\ell)\\
&=&\sum_{k\in G}\sum_{\ell\in G}f(k)f(\ell k)\rho(\ell k)\rho(k^{-1})\\
&=&\sum_{k\in G}f(k)\left(\sum_{\ell\in G}f(\ell k)\rho(\ell k)\right)\rho(k^{-1})\\
&=&\sum_{k\in G}f(k)\left(\sum_{s\in G}f(s)\rho(s)\right)\rho(k^{-1})\\
&=&\sum_{k\in G}f(k)\hat f(\rho)\rho(k^{-1})\\
&=&\hat f(\rho)\sum_{k\in G}f(k)\rho(k^{-1})\\
&=&\hat f(\rho)\sum_{k\in G}g(k^{-1})\rho(k^{-1})\\
&=&\hat f(\rho)\sum_{t\in G}g(t)\rho(t)\\
&=&\hat f(\rho)\hat g(\rho)
\end{eqnarray*}
\qed

Define the {\it temporal autocorrelation} by
\[ b_f(\ell)={\mathbb E}(f(v(T))f(v(T+\ell))) \]
for $\ell\in{\mathbb N}$ where the choice of $T\in{\mathbb N}$ is
immaterial since the random walk is stationary.

We can relate the spatial and temporal autocorelations by the following.
\begin{proposition}
\label{fourone}
\[ b_f(\ell)=\frac{1}{|G|^2}\sum_{\rho}d_{\rho} 
Tr(\hat\gamma(\rho)^{\ell}\hat h(\rho)) \]
where $h(s)=a_f(s^{-1})$ for $s\in G$ and the sum is over all irreducible
representations $\rho$ up to equivalence.
\end{proposition}

{\it Proof:}
\begin{eqnarray*}
b_f(\ell)&=&{\mathbb E}(f(v(T))v(T+\ell)))\\
&=&\frac{1}{|G|}\sum_{k\in G}{\mathbb E}(f(v(T))f(v(T+\ell))|v(T)=k)\\
&=&\frac{1}{|G|}\sum_{k\in G}\sum_{x\in G}f(k)\gamma^{*\ell}(x)f(xk)\\
&=&\frac{1}{|G|}\sum_{x\in G}a_f(x)\gamma^{*\ell}(x)
\end{eqnarray*}
where $\gamma^{*\ell}$ is the $\ell$-fold convolution of $\gamma$ with
itself, i.e. $\gamma^{*\ell}(k)={\mathbb P}(v(t+\ell)v(t)^{-1}=k)$. By the
Plancherel formula (as on p. 13 of Diaconis~\cite{di}),
\begin{eqnarray*}
\sum_{x\in G}a_f(x)\gamma^{*\ell}(x)&=&\sum_{x\in G}h(x^{-1})\gamma^{*\ell}(x)\\
&=&\frac{1}{|G|}
\sum_{\rho}d_{\rho}Tr({\widehat{\gamma^{*\ell}}}(\rho)\hat h(\rho))\\
&=&\frac{1}{|G|}\sum_{\rho}d_{\rho}Tr(\hat\gamma(\rho)^{\ell}\hat h(\rho))
\end{eqnarray*}
since 
${\widehat{\gamma^{*\ell}}}(\rho)=\hat\gamma(\rho)^{\ell}$. The proposition
follows.
\qed

Define the {\it multispectrum}
\[
A_f(\ell_1,...,\ell_n)=\sum_{k\in G}f(k)f(\ell_1 k)...f(\ell_n...\ell_1 k)
\]
for $\ell_1,...,\ell_n\in G$.

Define the {\it temporal multispectrum}
\[
B_f(\ell_1,...,\ell_n)={\mathbb E}(f(v(T))f(v(T+\ell_1))...f(v(T+\ell_1+...+\ell_n))).
\]

The Fourier transforms of $A_f$ and $B_f$ are related by the following.
\begin{proposition}
\label{fourmany}
\[
B_f(\ell_1,...,\ell_n)=\frac{1}{|G|^{n+1}}\sum (\prod_{i=1}^n d_{\rho_i})
Tr((\hat\gamma(\rho_1)^{\ell_1}\otimes ...\otimes \hat\gamma(\rho_n)^{\ell_n})
\hat H_n(\rho_1,...,\rho_n))
\]
where 
$H_n(x_1,...,x_n)=A_f(x_1^{-1},...,x_n^{-1})$,
$\hat H_n(\rho_1,...,\rho_n)$ is defined to be
$\hat H_n(\rho_1\otimes ...\otimes \rho_n)$, and the sum
 is such that $\rho_1,...,\rho_n$ each range over all irreducible
representations of $G$ up to equivalence.
\end{proposition}
Note that all irreducible representations of $G^n$ are, up to equivalence,
of the form $\rho_1\otimes ...\otimes\rho_n$ where $\rho_1,...,\rho_n$
are irreducible representations of $G$. (See, for example, p. 16 of 
Diaconis~\cite{di}.)

{\it Proof of Proposition~\ref{fourmany}:} 
Similarly to the proof of Proposition~\ref{fourone}, we get
\begin{eqnarray*}
\lefteqn{B_f(\ell_1,...,\ell_n)}\\
&=&
{\mathbb E}(f(v(T))f(v(T+\ell_1))...f(v(T+\ell_1+...+\ell_n)))\\
&=&\frac{1}{|G|}\sum_{k\in G}
{\mathbb E}(f(v(T))f(v(T+\ell_1)...f(v(T+\ell+1+...+\ell_n))|v(T)=k)\\
&=&\frac{1}{|G|}\sum_{k,x_1,...,x_n\in G}f(k)\gamma^{*\ell_1}(x_1)f(x_1k)...
\gamma^{*\ell_n}(x_n)f(x_n...x_1 k)\\
&=&\frac{1}{|G|}\sum_{x_1,...,x_n\in G}
\gamma^{*\ell_1}(x_1)...\gamma^{*\ell_n}(x_n)\sum_{k\in G}f(k)f(x_1k)...
f(x_n...x_1 k)\\
&=&\frac{1}{|G|}\sum_{x_1,...,x_n\in G} A_f(x_1,...,x_n)\gamma^{*\ell_1}(x_1)
...\gamma^{*\ell_n}(x_n)
\end{eqnarray*}

We shall use the Plancherel formula on $G^n$. First define
$p(x_1,...,x_n)=\gamma^{*\ell_1}(x_1)...\gamma^{*\ell_n}(x_n)$. Thus
\begin{eqnarray*}
B_f(\ell_1,...,\ell_n)&=&\frac{1}{|G|}\sum_{x_1,...,x_n\in G}H_n(x_1^{-1},
...,x_n^{-1})p(x_1,...,x_n)\\
&=&\frac{1}{|G|}\frac{1}{|G|^n}\sum_{\rho}d_{\rho}Tr(\hat H_n(\rho)\hat p(\rho))\\
&=&\frac{1}{|G|^{n+1}}\sum_{\rho}d_{\rho}Tr(\hat p(\rho)\hat H_n(\rho))
\end{eqnarray*}
where the sum is over all irreducible representations $\rho$ of $G^n$
up to equivalence. Such representations may be written in the form
$\rho=\rho_1\otimes ...\otimes\rho_n$. Then
\begin{eqnarray*}
\hat p(\rho)&=&\sum_{x_1,...,x_n\in G}\gamma^{*\ell_1}(x_1)...
\gamma^{*\ell_n}(x_n)\rho_1(x_1)\otimes ...\otimes\rho_n(x_n)\\
&=&\left(\sum_{x_1\in G}\gamma^{*\ell_1}(x_1)\rho_1(x_1)\right)\otimes ...
\otimes\left(\sum_{x_n\in G}\gamma^{*\ell_n}(x_n)\rho_n(x_n)\right)\\
&=&{\widehat{\gamma^{*\ell_1}}}(\rho_1)\otimes ...\otimes
{\widehat{\gamma^{*\ell_n}}}(\rho_n)\\
&=&\hat\gamma(\rho_1)^{\ell_1}\otimes ...\otimes\hat\gamma(\rho_n)^{\ell_n}
\end{eqnarray*}
The proposition follows. \qed

Linearity of the Fourier transform implies that to finish the
proof of Theorem~\ref{condition}, all we need to show is that $A_f$ suffices
to recover $f$ up to a shift, i.e. the following proposition.
\begin{proposition}
Suppose $A_{f_1}=A_{f_2}$ with $n=|G|$. Then $f_1$ is a shift of $f_2$.
\end{proposition}

{\it Proof:} First note that $A_f(x_1,...,x_n)>0$ if and only if there exists
an element $k\in G$ such that $f(k)=f(x_1k)=...=f(x_n...x_1 k)=1$. Number the
elements of $G$ from $1$ to $n$ such that the identity element $e$ is
numbered $n$. To an $n$-tuple $(x_1,...,x_n)\in G^n$, assign an $n$-tuple
$(a_1,...,a_n)$ of integers so that $a_1$ is the number of $x_1$ and
if $2\le j\le n$, $a_j$ is the smallest integer whuch is greater than
$a_{j-1}$ and congruent modulo $n$ to the number of $x_j...x_1$. Let
$m(f)=(m_1(f),...,m_n(f))$ satisfy $A_f(m_1(f),...,m_n(f))>0$ such that
the $n$-tuple $(a_1,...,a_n)$ assigned to it is the lexicographically smallest
$n$-tuple assigned to an $n$-tuple $(x_1,...,x_n)\in G^n$ with
$A_f(x_1,...,x_n)>0$. (If there are no $n$-tuples $(x_1,...,x_n)$ with
$A_f(x_1,...,x_n)>0$, then $A_f(e,...,e)=0$ where $e$ is the identity element
of $G$ and $f(k)f(ek)...f(e^nk)=0$ and hence $f(k)=0$ for all $k\in G$.)

Let $i$ be the largest index such that $a_i<n$ where $(a_1,...,a_n)$ is
assigned to $(m_1(f),...,m_n(f))$. For some $k\in G$, $f(k)$, $f(m_1k)$,
..., $f(m_i(f)...m_1(f)k)$ are all $1$; otherwise $A_f(m_1(f),...,m_n(f))$
would be $0$. Now suppose $f(x)=1$ for some 
$x\notin\{k,m_1(f)k,...,m_i(f)...m_1(f)k\}$. Let $y$ be such that $x=yk$, i.e.
$y=xk^{-1}$. We shall create an $n$-tuple $K=(k_1,...,k_n)$ of elements of $G$
such that $A_f(k_1,...,k_n)>0$ while the $n$-tuple of integers assigned to
$K$ is lexicographically smaller than the $n$-tuple of integers assigned to
$m(f)$, contradicting the definition of $m(f)$. If the number assigned to $y$
is greater than the number assigned to $m_i(f)...m_1(f)$, then let
$k_1=m_1(f)$, ..., $k_i=m_i(f)$, $k_{i+1}=y(m_i(f)...m_1(f))^{-1}$, $k_{i+2}=e$,..., $k_n=e$ where $e$ is the identity element of $G$. Otherwise let $j$ be the smallest value such that the number
assigned to $y$ is less than the number assigned to $m_j(f)...m_1(f)$. Let
$k_1=m_1(f)$, ..., $k_{j-1}=m_{j-1}(f)$, $k_j=y(m_{j-1}(f)...m_1(f))^{-1}$,
$k_{j+1}=m_j(f)m_{j-1}(f)...m_1(f)y^{-1}$, $k_{j+2}=m_{j+1}(f)$, ..., 
$k_{i+1}=m_i(f)$, $k_{i+2}=e$, ..., $k_n=e$. (In particular, if $j+2\le b\le
i+1$, then $k_b=m_{b-1}(f)$.) In either case, it can be verified that
$A_f(k_1,...,k_n)>0$ while the $n$-tuple of integers assigned to $K$ is 
lexicographically smaller than the $n$-tuple of integers assigned to $m(f)$.
This contradiction implies that $f(x)=0$ if $x$ is not one of $k$, $m_1(f)k$,
..., $m_i(f)...m_1(f)k$, and so $A_f$ determines $f$ up to a shift.

The proposition follows, and so does Theorem~\ref{condition}.\qed

\section{Proof of Theorem~\ref{muchado}}

Each irreducible representation is equivalent to an irreducible
representation $\rho$ such that $\hat\gamma(\rho)$ is upper triangular (and
in Jordan canonical form). Thus we may without loss of generality assume that
$\hat\gamma(\rho_1)$, ..., $\hat\gamma(\rho_n)$ are all upper triangular. The
elements of $\hat J_n(\rho_1\otimes...\otimes\rho_n)$ are linear combinations
of $J_n(x_1,...,x_n)$ where $(x_1,...,x_n)$ range over the $n^n$ elements of 
$G^n$. If $\hat\gamma(\rho_1)$,..., $\hat\gamma(\rho_n)$, and hence (with
a natural basis) $\hat\gamma(\rho_1)^{\ell_1}\otimes...\otimes
\hat\gamma(\rho_n)^{\ell_n}$ are upper triangular, then
$Tr((\hat\gamma(\rho_1)^{\ell_1}\otimes...\otimes\hat\gamma(\rho_n)^{\ell_n})
\hat J_n(\rho_1\otimes...\otimes\rho_n))$ excludes elements above the diagonal
of $\hat J_n(\rho_1\otimes...\otimes\rho_n)$. If $d_{\rho_i}>1$ for some $i$,
then there will be elements above the diagonal. If $G$ is a finite non-abelian
group, then $d_{\rho_i}>1$ for some irreducible representation $\rho_i$. (See,
for example, p. 15 of Diaconis~\cite{di}.) Also $\sum_{\rho_i}d_{\rho_i}^2=n$
where the sum is over all irreducible representations of $G$ up to equivalence. The total number of elements for all matrices $\hat J_n(\rho_1\otimes ...\otimes\rho_n)$ is $\sum_{\rho_1}...\sum_{\rho_n}(d_{\rho_1}...d_{\rho_n})^2=n^n$.
When $G$ is a non-abelian group, the equations
\[
\sum d_{\rho_1}...d_{\rho_n}Tr((\hat\gamma(\rho_1)^{\ell_1}\otimes...\otimes
\hat\gamma(\rho_n)^{\ell_n})\hat J_n(\rho_1\otimes...\otimes\rho_n))=0
\]
over all positive integers $\ell_1,...,\ell_n$ 
give rise to a system of homogeneous linear equations involving only
elements which are on or below the diagonal of $\hat J_n(\rho_1\otimes...\otimes\rho_n)$
for some $\rho_1,...,\rho_n$. When all the elements which are on or below this diagonal
for some $\rho_1,...,\rho_n$ are $0$, this system of equations is satisfied.
Since there are less than $n^n$ such elements if $G$ is non-abelian, solutions
exist where not all $J_n(x_1,...,x_n)$ are $0$, and Theorem~\ref{muchado}
follows. \qed

\section{Questions for Further Study}

In addition to other questions posed in \cite{fty}, the work here leaves open
the question if there are reconstructive random walks on finite non-abelian
groups. Perhaps the solutions where $J_n$ is not identically $0$ do not come 
from the difference of two multispectrums $A_{f_1}$ and $A_{f_2}$ of sceneries.
Computer exploration on small non-abelian groups might be a place to start
exploring that question.
Indeed some computer exploration with Maple encouraged the author to consider
Theorem~\ref{muchado}.


\begin{thebibliography}{99}

\bibitem{di} P. Diaconis. {\it Group Representations in Probability and
Statistics}. Hayward, Calif: Institute of Mathematical Statistics, 1988.

\bibitem{fty} H. Finucane, O. Tamuz, and Y. Yaari. 
Scenery reconstruction on finite abelian groups. Stochastic Process. Appl.
124 (2014) 2754-2770.

\bibitem{ml} H. Matzinger and J. Lember. Reconstruction of periodic sceneries
seen along a random walk. Stochastic Process. Appl. 116 (2006) 1584-1599.

\bibitem{se} J.P. Serre. {\it Linear Representations of Finite Groups}.
New York: Springer-Verlag, 1977.

\end{thebibliography}
\end{document}